\documentclass[11pt]{article}
\usepackage{mathrsfs}
\usepackage{amsfonts}
\usepackage{bbm}
\usepackage{latexsym,amsfonts,euscript,amssymb}

\newtheorem{lem}{\bf Lemma}[section]

\newtheorem{defi}[lem]{\bf Definition}

{\left\lbrace\matrix}{endmatix\right\rbrace}

 {\title{
 A Proof Of The Riemann Hypothesis
\thanks{Project supported by the National Natural Science Foundation of China(Grant No. 11871224).
E-mail:xumch@scnu.edu.cn
}}

\author{ Mingchun Xu\\  \small School of Mathematics,
  South-China Normal University, \\ Guangzhou, 510631, China }
\begin{document}
\maketitle

\begin{abstract}
   We consider the alternating
Riemann zeta function $\zeta^*(s)= \sum^{\infty} _{ n=1}
\frac{(-1)^{n-1}}{n^s}$,  which converges if $Re (s)>0 .$  By using  Rouche's theorem, the Bolzano-Weierstrass theorem  and by method  of contradiction  we complete
the proof of the Riemann Hypothesis.

\end{abstract}


\textbf{2010 MR Subject Classification}\,\, 11M06,11M41

\section{Introduction}

The distribution of prime numbers is an old problem in number
theory. It is very easy to state but extremely hard to resolve. In
his famous paper written in 1859 Bernhard Riemann connected this
problem with a function investigated earlier by Leonhard Euler. He
also formulated certain hypothesis concerning the distribution of
complex zeros of this function. At first this hypothesis appeared as
a relatively simple analytical conjecture to be proved sooner rather
than later. However, future development of the theory proved
otherwise: since then the Riemann hypothesis (hereafter called RH)
is commonly regarded as both the most challenging and the most
difficult task in number theory. It states that all complex zeros of
the zeta function, defined by the following series if  Re $ s > 1$
$$\zeta(s) = \sum^{\infty} _{ n=1} \frac{1}{n^s}
$$
 and by analytic continuation to the whole plane, are located right on the critical line
$Re(s) = \frac{1}{2}.$  %For a  rich history of the Riemann
%hypothesis and some recent developments, see Bombieri\cite{bo},
%Conrey\cite{c}, and Sarnark\cite{s}.
We denote by $\mathbb{C}$,  $Re(s)$ and $Im(s)$ the complex  $s$-plane, the real and
the imaginary part of the complex number $s$ respectively.
We have already seen that the zeros of $\zeta(s)$ at $-2, -4, \cdots $, are known as the trivial zeros and the nontrivial zeros are all complex in the strip $\{s| 0<Re(s)<1 \}$.
Hence the nontrivial zeros either lie on   the critical line
$Re(s) = \frac{1}{2} $, or occur in pairs symmetrical about this line and the real axis by the functional equation.
In this paper,  we prove the following result.

 \textbf{Main  Theorem }\,\,\  All nontrivial
zeros of the  Riemann zeta function $\zeta(s)$  lie  on the critical line
$Re(s) = \frac{1}{2} .$

\section{  %Compactness and Convergence in the Space of Analytic Functions
Some definitions and lemmas}

In this paper a metric is put on the set of all analytic functions
on a fixed region $G$, and compactness and convergence in this
metric space is discussed. For further detail see\cite  {con}.
Let $(\Omega, d)$ denote a complete metric space.
If $G$ is an open set in $\mathbb{C}$ and
$(\Omega, d)$ is  a complete metric space then $ C(G, \Omega)$
denotes the set all continuous functions from $G$ to $\Omega$.

\begin{lem} (see\cite{con},Chapter VII,1.2 Proposition) If $G$ is an open set in $\mathbb{C}$ then there is sequence $\{K_n\}$ of
compact subsets of  and $G$ such that $G=\bigsqcup^\infty_{n=1}
K_n$. Moreover, the sets $K_n$ can be chosen to satisfy the
following conditions:

(a)\,\, $K_n \subset int\,\ K_{n+1}$;

(b)\,\, $K\subset G$ and $K$ compact implies $K\subset K_n$ for some
$n$;

(c)\,\, Every component of $\mathbb{C}_\infty-K_n$ contains  a
component of $\mathbb{C}_\infty-G.$
\end{lem}

\textbf{Definition} If $G=\bigsqcup^\infty_{n=1} K_n$ where
each$\{K_n\}$  is  compact  and $K_n \subset int\,\ K_{n+1}$ as the
above lemma, define    $$\rho_n(f,g)=sup\{ d(f(z),g(z))| z\in
K_n\}$$ for all functions $f$ and $g$ in $C(G, \Omega)$. Also define
$$\rho(f,g)=\sum_{n=1}^\infty (\frac{1}{2})^n \frac{\rho_n(f,g)}{1+\rho_n(f,g)}.$$
\begin{lem} (see\cite{con},Chapter VII,1.6 Proposition)
 $(C(G, \Omega), \rho)$ is metric space.
\end{lem}

\begin{lem}(see\cite{con},Chapter VII,1.10 Proposition)
A sequence $\{ f_n \}$ in $(C(G, \Omega), \rho)$ converges to $f$ if
and only if  $\{ f_n \}$ converges to $f$ uniformly on all compact
subsets of $G$.
\end{lem}
\textbf{Definition} (see\cite{R},Chapter 1, \S 1.1)For the above lemma the sequence $\{ f_n \}$ in \linebreak $(C(G, \Omega), \rho)$, which converges to $f$ with respect to the metric $\rho$,  is called compactly convergent in  $G$.

%\begin{defi}
  If $G$ is a region in $\mathbb{C}$ then
$H(G)$ denotes the set all analytic functions. We will always assume
that the metric on $H(G)$ is the metric which it inherits as a
subset of $C(G,\Omega).$
%\end{defi}

\begin{lem}(see\cite{R},Chapter 7, \S 1.2 footnote or see\cite{con},Chapter VII, \S 1, exercise 5)

 Let $f_n(s)$ be a sequence of functions, each analytic  in a region $G$,
 and let $f_n(s)$ be compactly convergent to $f(s)$ in  $H(G)$.  If  points $z_n$ has a limit-point $z$
inside $G$, then $\lim_{n\rightarrow \infty}f_n(z_n)=f(z)$.
\end{lem}
\begin{defi} A set $\mathscr{F}\subset H(G)$ is locally bounded if for each point $a$ in $G$
there are constant $M$ and $r>0$ such that for all $f$ in
$\mathscr{F}$,  $$|f(z)|\leq M, \mbox{for}\, \, |z-a|<r. $$

Alternately,  $\mathscr{F}\subset H(G)$ is locally bounded if there
is an $r>0$ such that
$$sup\{ |f(z)|:|z-a|<r, f\in \mathscr{F} \}<\infty.
$$ That is,  $\mathscr{F}$ is locally bounded if about  each point $a$ in $G$
there is disk on which $\mathscr{F}$ is uniformly bounded. This
immediately extends to the requirement that $\mathscr{F}$ be
uniformly bounded on compact sets in $G$.
\end{defi}
\begin{lem} A set $\mathscr{F}\subset H(G)$ is locally bounded iff for each compact set  $K\subset G$ there is
 constant $M$  such that for all $f$ in $\mathscr{F}$ and $z$ in $K$
$$|f(z)|\leq M.$$ \end{lem}
%\begin{thm} (Montel's Theorem)\,\

%A family  $\mathscr{F}\subset H(G)$ is normal iff $\mathscr{F}$ is
%locally bounded. \end{thm}
\begin{lem} (Vitali's convergence theorem, see\cite{T1}, 5.21)\,\,
 Let $f_n(s)$ be a sequence of functions, each analytic  in a region $G$;
 let  $\{f_n(s)\}\subset H(G)$ is locally bounded;
 and let $f_n(s)$ tend to a limit, as
$n\longrightarrow \infty, $ at a set of points having a limit-point
inside $G$. Then $f_n(s)$ tends uniformly to a limit in any region
bounded by a contour interior to $G$, the limit being, therefore, an
analytic function of $s$.
\end{lem}

\begin{lem} ( see\cite{con},Chapter V, 3.8 Rouch$\acute{e}$'s Theorem)\,\ Let $\gamma$ denote  a simple closed contour, and suppose that

(a)\,\,\  two functions  $f(s)$ and $g(s)$ are  analytic inside and on  $\gamma$;

(b)\,\,  $|f(s)|>|g(s)|$ on each point on  $\gamma$.

Then $f(s)$ and  $f(s)+g(s)$ have  the same number of zeros, counting multiplicities, inside  $\gamma$.  \end{lem}

\begin{lem} (Hurwitz's theorem, see\cite{R},Chapter 7, \S 5)\,\,
 Let $G$ be  a region $G$ in $\mathbb{C}$, and let $f_n(s)$ be a sequence of functions, each analytic  in $G$, that  converges compactly to $f(s)\in  H(G)$. If all the  functions $f_n(s)$ are nonvanishing
 in $G$ and $f(s)$ is not identically zero, then  $f(s)$ is nonvanishing in $G$.
\end{lem}

%\begin{lem} ( see\cite{T2},Chapter 3, \S 3.15)\,\, Let  $\mathcal {P}$ be the set of all prime numbers.   If $z=1+t\, i$ and $0\neq t\in \mathbf{R}$,  then $$\sum_{p\in \mathcal {P}}p^{-z}$$ is convergent.
%\end{lem}
\section{Proof of the Main Theorem }

The Riemann zeta
function $\zeta(s)$ is analytic in  the complex plane ($s$-plane),
except for a simple pole at $s=1.$ The alternating zeta function
$\zeta^*(s)$ is defined as the analytic continuation of the
Dirichlet series
\begin{equation}
\label{f-1} \zeta^*(s) = \sum^{\infty} _{ n=1}
\frac{(-1)^{n-1}}{n^s},
 \end{equation}
  which converges if $Re(s)>0 .$
The two functions are related to for all complex $s$ by the identity
\begin{equation}
\label{f-1}
 \zeta^*(s)=(1-2^{1-s})\zeta(s),
 \end{equation}
  which is easily established,
first for $Re(s)>1$ by combining terms in the convergent
Dirichlet series, and then by using analytic continuation to extend
the result $s$-plane. The factor $1-2^{1-s}$ has a simple zero at
$s=1$ that cancels the pole of $\zeta(s)$, so (2) shows that
$\zeta^*(s)$ is an entire function of $s$. It vanishes at each zero
of the factor $1-2^{1-s}$ with the exception of $s=1$, at which
point $\zeta^*(1)=ln\, 2.$ The other zeros of factor $1-2^{1-s}$ are
called also trivial zeros of $\zeta^*(s)$, which lie on the  line Re
$s =1$ and occur at the points at $s=1+\frac{2l\pi}{ln 2} i$ for all
nonzero integers $l$. All nontrivial zeros of $\zeta(s)$ are the
same as that of $\zeta^*(s)$.

In this section we  will use  Rouch$\acute{e}$'s theorem  to prove the Riemann Hypothesis. The proof will be by contradiction.
Since $ \zeta^*(s)= \sum^{\infty} _{ n=1}
\frac{(-1)^{n-1}}{n^s} >0\,\, ( s\in \mathbb{R}, s>0),$
$\zeta^{*}(s)$ has no zeros on the positive axis.
 Assume there is a zero $b$ of $\zeta^{*}(s)$
 in the interior of the region $\{ s|\frac{1}{2}< Re(s)<1, Im(s)> 0 \} $. We can similarly deal with the case of a zero $b$ of $\zeta^{*}(s)$
 in the interior of the region $\{ s|\frac{1}{2}< Re(s)<1, Im(s)< 0 \} $.  To prove Riemann Hypothesis, it is then sufficient to show that the existence of $b$ leads a contradiction. Then there is some compact set $F\subset \{ s|\frac{1}{2}< Re(s), Im(s)> 0 \} $ with a  simple closed
contour $\gamma$ as the boundary of $F$ such that $b$ is in the interior of $F$ and $b$ is just only one zero of  $\zeta^{*}(s)$, where $\gamma=\gamma_{r}\cup\gamma_{l}$, $\gamma_{l}$
 is the left part of $\gamma$ in the set  $\{ s|\frac{1}{2}< Re(s)\leq 1, Im(s)> 0 \} $ and $\gamma_{r}$
 is the right part of $\gamma$ in the set  $\{ s| Re(s)\geq 1, Im(s)> 0 \} $, $\gamma \cap \{ s| Re(s)> 1, Im(s)> 0 \} \neq \varnothing $.
 The line $Re(s)=1$ intersects with $\gamma$ at the points $A, B$ respectively such that $A, B$ are between two points  $1+ \frac{2l\pi}{ln 2} i$ and  $1+ \frac{2(l+1)\pi}{ln 2} i$ for some positive integer $l$ and $\{ A,B\} \cap \{1+ \frac{2l\pi}{ln 2} i,1+ \frac{2(l+1)\pi}{ln 2} i\} =\emptyset $ ( in FIGURE 1).

 \setlength{\unitlength}{1.5cm}
\begin{picture}(5,3)
\put(.5,0.2){\vector(1,0){8}}
\put(6,0.5){\circle*{.1}}\put(6.5,.35){ 1+ $\frac{2l\pi}{ln 2} i$}\put(6,0.5){\line(1,0){.9}}
 \put(6,2){\circle*{.1}}\put(6,2.5){ 1+ $\frac{2(l+1)\pi}{ln 2} i$}\put(6,2){\line(1,1){0.5}}

\put(4.5,-0.1){\line(0,1){3}}
\put(6,-0.1){\line(0,1){3}}
\put(3,-0.1){\vector(0,1){3}}

%\put(5.3,1.1){\line(0,1){0.35}}
%\qbezier(5,1)(6.1,1.7)(7,1)
%\qbezier(5,1)(6.1,0.7)(7,1)
%\qbezier(4.8,1)(6.1,1.9)(7.2,1)
%\qbezier(4.8,1)(6.1,0.5)(7.2,1)
\qbezier(4.7,1)(6.1,2.2)(7.3,1)
\qbezier(4.7,1)(6.1,0.3)(7.3,1)

\put(6,0){\parbox{0.5cm}{1}}
\put(4.5,0)
{\parbox{0.5cm}{0.5}}
\put(3,0.1){\parbox{0.5cm}{O}}

\put(7.6,1.3){\parbox{0.5cm}{$\gamma$ }}
\put(7.6,1.3){\line(-1,-1){0.4}}
%\put(7.3,1.7){\parbox{0.5cm}{$\gamma_t$ }}\put(7.3,1.7){\line(-1,-1){0.45}}
%\put(4.5,.3){\parbox{0.5cm}{$\gamma_1$ }}\put(4.8,.45){\line(1,1){0.46}}
\put(6,.55){\parbox{0.5cm}{$A$}}\put(6,.65){\circle*{.1}}
%\put(6,.75){\circle*{.1}}\put(6,.85){\circle*{.1}}\put(6,.85){\parbox{0.5cm}{$A_{1}$}}\put(5.7,.75){\parbox{0.5cm}{$A_{t}$}}

\put(6,1.7){\parbox{0.3cm}{$B$}}
%\put(6,1.6){\line(1,1){0.3}}
\put(6,1.6){\circle*{.1}}

%\put(6,1.45){\circle*{.1}}\put(6.5,1.6){\parbox{0.3cm}{$B_{t}$}}\put(6,1.45){\line(2,1){0.3}}
%\put(6,1.35){\circle*{.1}}\put(6,1.2){\parbox{0.3cm}{$B_{1}$}}

%\put(4.5,2){\parbox{0.3cm}{$D_{0}$}}
%\put(4,1.8){\parbox{0.3cm}{$D_{t}$}}\put(3.8,1.4){\parbox{0.3cm}{$D_{1}$}}
%\put(5.3,1.4){\line(-1,1){0.5}}
%\put(5.3,1.45){\circle*{.1}}\put(5.3,1.28){\circle*{.1}}\put(5.3,1.28){\line(-2,1){1}}\put(5.3,1.13){\line(-3,1){1}}
\put(5.4,1.6){\parbox{0.3cm}{$C$}}
%\put(5.3,1.13){\circle*{.1}}
\put(5.4,1.45){\circle*{.1}}
%\put(5.45,1.35){\circle*{.1}}\put(5.46,1.25){\circle*{.1}}
%\put(5.43,1.23){\line(0,1){0.25}}

\put(5.6,0.9){\parbox{0.5cm}{$b$}}\put(5.5,1){\circle*{.1}}
%\put(5.85,1.45){\line(0,1){1.5}}\put(5.85,3){\parbox{0.5cm}{$\Gamma$ }}
\put(0,1){\parbox{1cm}{$FIGURE \, \, 1$}}
\end{picture} }

Denote by $\mathcal {P}$, $\mathbb{N}$  the set of all prime numbers and the set of all positive integers respectively. Let $n\in \mathbb{N}$ and subdivide a sequence 
\begin{equation}
g_n(s)=(\prod_{p\leq n, p\in \mathcal {P}} exp (-p^{-s}) ) \cdot\zeta^{*}(s) = f_n(s)+ r_n(s)
\end{equation}
into the main term
 $$ f_n(s)=(1-2^{1-s})\prod_{p\leq n, p\in \mathcal {P}}(1-p^{-s})^{-1} exp (-p^{-s})$$ plus the remainder term $r_n(s)$.

 Step 1\,\, For $Re(s)>\frac{1}{2}$ we have $$\eta(s) =(1-2^{1-s})\prod_{ p\in \mathcal {P}}(1-p^{-s})^{-1} exp (-p^{-s})$$ is
analytic.

% If $Re(s)\geq a >\frac{1}{2}$, then
%$$| ln((1-p^{-s})exp(p^{-s})) |=|p^{-s}+ln(1-p^{-s})|=|p^{-s}-p^{-s}-\frac{1}{2}p^{-2s}-\frac{1}{3}p^{-3s}-\cdots|$$
%$$=|-\frac{1}{2}p^{-2s}-\frac{1}{3}p^{-3s}-\cdots|\leq \frac{1}{2}(|p^{-2s}|+|p^{-3s}|+\cdots)=\frac{1}{2}\frac{|p^{-2s|}}{1-|p^{-s}|}\leq \frac{1}{p^{2a}}$$
%for all sufficiently large values of $p$ and with
%$|p^{-s}|<\frac{1}{2}.$  Note that $\sum_{p\in\mathcal {P} }p^{-2a}$
%is convergent for $ a >\frac{1}{2}$.

As $u\rightarrow 0$, $(1-u)^{-1}e^{-u}=1+O(u^2)$. Since the series $\sum_{p\in \mathcal {P}}
|p^{-2s}|$ converges uniformly to a bounded sum, by $|p^{-2s}|\leq p^{-2a}, {Re}(s)= a, 2a>1$  and Weierstrass M-test, we have
$$\prod_{p\in\mathcal {P}}((1-p^{-s})^{-1}exp(-p^{-s}))$$ is uniformly and
absolutely convergent in region $ {Re}(s)\geq c
>\frac{1}{2}$ by \cite  {con}, Chapter 7,\S 5, theorem 5.9 or \S 1.43 and \S 1.44 in \cite{T1}.
 Then $f_n(s)$ is compactly
convergent to $\eta(s) $ in $H(\{s| Re(s)>\frac{1}{2}\})$.  Thus $\eta(s)$ is nonvanishing on  $\gamma$ by Lemma 2.9.

 Step 2\,\,
 Since  $f_n(s)$ and $g_n(s)$ converge uniformly to  $\eta(s) $ for each compact subset in the set $\{ s|Re(s)\geq 1, Im(s)> 0\}$  by $\cite{T2}$,Chapter 3, $\S 3.15$, we have
$r_n(s)$ converges uniformly to  $0 $ for each compact subset in the set $\{ s|Re(s) \geq 1, Im(s)> 0\}$.
%Thus we have
%\begin{equation}
%\lim_{n\rightarrow \infty}|r_n(s_{00})|=\lim_{n\rightarrow \infty}|r_n(s_{10})|= 0.
%\end{equation}

Thus there is a $N\in  \mathbb{N}$ such that if $n> N$ then $|f_n(s)|>|r_n(s)|$ holds uniformly  for each $s\in \gamma_r=\gamma\cap \{ z|Re(z)\geq 1, Im(z)> 0\}$.

% Step 3\,\,
%  Since $b$ is just only one zero of  $\zeta^{*}(s)$   inside to $\gamma$,  and $|f_n(s)|-|r_n(s)|$ is continuous in  $\{ s| Re(s)> \frac{1}{2}, Im(s)> 0 \} $,  by Rouch$\acute{e}$'s theorem and intermediate value theorem of continuous function,  for  $n> N$ there are  $z_n$ such that
 %\begin{equation}|f_n(z_n)|=|r_n(z_n)| \end{equation}
%   holds but
% $|f_n(s)|>|r_n(s)|$ holds uniformly  for each $s\in \gamma_{l}$ and $s$ being between $z_n$ and $B$.
 % Now we note that $\gamma_{l}$ is a bounded closed set and the sequence  $Z=\{z_n|n>N\}$   is infinite. Then $Z$  has  a limit point
% $C \in \gamma_{l}$ by the Bolzano-Weierstrass theorem. Thus there is a subsequence $n_k$ of  $\mathbb{N}$ such that $$  \lim_{k\rightarrow \infty}z_{n_k}=C.$$
% Thus $\{g_{n_k}(z_{n_k})|n_k>N\}$ is a bounded sequence  by (3) and (4) and Step 1.  By the Bolzano-Weierstrass theorem  there is a subsequence $m_k$ of  $n_k$ such that $$  \lim_{k\rightarrow \infty}g_{m_k}(z_{m_k})$$ exists.
 %Thus  $$\sum_{p\leq m_k,p\in \mathcal {P}}
%p^{-z_{m_k}}$$ converges by (3) and $\zeta^{*}(C)\neq 0$. If $C$ is  either $B$ or $A $,  then $$\sum_{p\in \mathbb{P}}p^{-C}$$ is convergent by Lemma 2.7 and by (4)  we have
% \begin{equation} \eta(C)= \lim_{k\rightarrow \infty}r_{m_k}(z_{m_k})= 0,\end{equation} which is a
% contradiction to Step 1.

 Step 3\,\, Since $b$ is just only one zero of  $\zeta^{*}(s)$   in the interior of $F$,  and $|f_n(s)|-|r_n(s)|$ is continuous in  $\{ s|\frac{1}{2}< Re(s), Im(s)> 0 \} $,  by Rouch$\acute{e}$'s theorem and intermediate value theorem of continuous function,  for  $n> N$ there is a $z_n\in  \gamma_{\emph{l}}$ such that
 \begin{equation}|f_n(z_n)|=|r_n(z_n)|\end{equation}
   holds but
 $|f_n(s)|>|r_n(s)|$ holds uniformly  for each $s\in \gamma_{\emph{l}}$ and $s$ being between $z_n$ and $B$.
  Now we note that $\gamma_1$ is a bounded closed set and the sequence  $Z=\{z_n|n>N\}$ is infinite. Then $Z$ has a limit point
 $C\in \gamma_1$. Thus there is a subsequence $z_{n_k}$ of $Z$ such that $$ \lim_{k\rightarrow \infty}z_{n_k}=C .$$

 By Lemma 2.4 and (4) we have
 \begin{equation} \lim_{k\rightarrow \infty}|r_{n_k}(z_{n_k})|=|\eta(C)|\neq 0.\end{equation}
 %By (4) and (6) we get $s_3\neq s_1$ and   $s_3\neq s_2$.
 Furthermore, if $n_k>N$,  then $g_{n_k}(s)=f_{n_k}(s)+O(f_{n_k}(s))$ holds uniformly  for each $s\in \gamma_1$ and $s$ being  between $C$ and $B$.

Step 4 \,\,% Next  $C$ in  $\gamma _l$ is neither $B$ nor $A $.
 Denote by $\gamma_3$ the bounded closed set $\{s|s\in \gamma_1$ and $s$ is between $C$ and $B$ including $C$ and $B \}$.
From Step 3 for  $n_k> N$ there is a $B(s;r)=\{z| |z-s|<r\} \subset \{ s|Re(s)>\frac{1}{2}, Im(s)> 0\}$, for some $r>0$ , $s\in\gamma_3$ such that $g_{n_k}(z)=f_{n_k}(z)+O(f_{n_k}(z))$  holds uniformly  for each $z\in B(s;r) $.

By the  Heine-Borel finite covering theorem  there are $m$ points \linebreak $\{ a_1,a_2, \cdots, a_m\}\subseteq \gamma_3$ and  $B(a_i;r_i)=\{z| |z-a_i|<r_i\}$ for some $r_i>0$ , $i=1,2,\cdots, m$ such that  $g_{n_k}(z)=f_{n_k}(z)+O(f_{n_k}(z))$  holds uniformly  for each $z\in B(a_i;r_i) $,$i=1,2,\cdots, m$, $n_k> N$ and $\bigcup_{i=1}^{m}B(a_i;r_i)\supseteq \gamma_3 $.

 Let $G=\bigcup_{i=1}^{m}B(a_i;r_i)\cup \{ s|Re(s)>1,Im(s)>0 \}$. Since $f_n(s)$ is
convergent to $\eta(s) $ in $H(\{s| Re(s)>\frac{1}{2}\})$, we have $g_{n_k}(s)$ are locally bounded in $G$ for $n_k>N$.

Since  $g_n(s)$ is
convergent to $\eta(s) $ in $H(\{ s|Re(s)>1,Im(s)>0 \})$, by Vatali's theorem, we have $g_{n_k}(s)$ is
convergent to $\eta(s) $ in $H(G)$.

Step 5 \,\, Since $f_{n_k}(s)$ and  $g_{n_k}(s)$ are
convergent to $\eta(s) $ in $H(G)$, by Lemma 2.4 and (3) we have
 \begin{equation} \lim_{k\rightarrow \infty}r_{n_k}(z_{n_k})=0.\end{equation}
 Finally we have (6) contradicts (5). This completes the proof of the theorem.

%\textbf{Acknowledgement} The author would like to thank Professor
%Wujie Shi for his guidance and  his valuable help. \noindent

\end{document}